\documentstyle[12pt]{article}
\title{\bf Oscillation of
a Linear Delay Impulsive Differential Equation}

\author{L. Berezansky, $^{1}$
\\Ben-Gurion University of the Negev,\\
Department of Mathematics and Computer Science,\\
Beer-Sheva 84105, Israel, \\
E. Braverman \\
Technion - Israel Institute of Technology, \\
Israel Institute of Metals, 32000, Haifa, Israel }

\begin{document}
\maketitle

\footnotetext[1]{Supported by the Israel Ministry
of Science and Technology and Israel Ministry of Absorption}

\begin{abstract}
The main result of the paper is that the oscillation (non-oscillation)
of the impulsive delay differential equation
$$
\dot {x}(t)+\sum _{k=1}^m A_k(t)x[h_k(t)]=0,~~t\geq 0,
$$
$$
x(\tau_j)=B_jx(\tau_j-0),~~\lim \tau_j = \infty
$$
is equivalent to the oscillation (non-oscillation)
of the equation without impulses
$$
\dot {x}(t)=\sum_{k=1}^m A_k(t)~\prod_{h_k(t)<\tau_j\leq t}
B_j^{-1}x[h_k(t)]=0,~~t\geq0.
$$

Explicit oscillation results
are presented.
\end{abstract}

\section{Introduction}

Recently results on
oscillation of delay differential
equations have taken shape of a developed theory
presented in monographs
[1-4].
At the same time it is an intensively developing field
which is an objective of numerous publications.

However, for impulsive differential equations
there are only few publications dealing with
oscillation problems [1,4,5,6].

The purpose of the present paper is to fill up this gap.
The main result is that the oscillation (non-oscillation)
of the impulsive delay differential equation
is equivalent to the oscillation (non-oscillation)
of a certain differential equation without impulses
which can be constructed explicitly from an impulsive equation.
Thus the oscillation problems (in particular, oscillation
and non-oscillation criteria) for an impulsive equation can be
reduced to the similar problem for a certain non-impulsive
equation.

The method proposed in the present paper for oscillation
is new both for impulsive and non-impulsive equations.
It is based on the solution representation formula.
Recently such formulas are widely used in stability
investigations of non-impulsive [7-9]
and impulsive equations [5,10-12].

We demonstrate that the existence of a nonoscillating
solution is equivalent to the positiveness of the fundamental
function.
At the same time this is equivalent to the solvability
of a certain nonlinear inequality
which is similar to "the generalized characteristic
equation" from the monograph [2].

The paper is organized as follows.
Theorems 1 and 2 are concerned with the equivalence
of non-oscillation, positiveness of a fundamental function
and solvability of a certain inequality.
They lead to explicit non-oscillation results (Theorem 3).
Theorem 4 compares non-oscillation conditions
for two different impulsive delay
differential equations.
Theorems 5 and 6 give new oscillation criteria
for delay differential equations without impulses.
Theorem 7 contains the main result of the paper
connecting oscillation of an impulsive and
a non-impulsive equation.
As a corollary (Theorem 8) we obtain explicit
oscillation conditions for an impulsive delay equation.

We are very grateful to Prof. Yury Domshlak
for useful discussions of the problems considered in the paper.

\section{Preliminaries}

We consider a scalar delay differential equation
\begin{equation}
\dot{x}(t) + \sum_{k=1}^m A_k(t) x[h_k(t)]=f(t),
{}~~~t \geq 0;
\end{equation}
\begin{equation}
x(\tau_j)=B_j x(\tau_j -0), ~j=1,2, \dots ,
\end{equation}
under the following assumptions

(a1) $0 = \tau_0 < \tau_1 < \tau_2 < \dots $
are fixed points, $\lim_{j \rightarrow \infty}\tau_j = \infty$;

(a2) $A_k, f, k=1, \dots, m$ are Lebesgue measurable
functions essentially bounded in each finite interval
$[0,b], ~B_j \in {\bf R}, ~j=1, \dots, ~{\bf R}$
is a real axis;

(a3) $h_k : [0, \infty) \rightarrow {\bf R}$
are Lebesgue measurable functions, $h_k(t) \leq t$.

Together with (1),(2) we will consider
for each $t_0 \geq 0$ an initial value problem
\begin{equation}
\dot{x} (t) +\sum_{k=1}^m A_k (t) x[h_k(t)] =f(t),
\mbox{~where~}t\geq t_0,
{}~x(\xi)=\varphi(\xi), \xi < t_0,
\end{equation}
\begin{equation}
x(\tau_j) = B_j x(\tau_j -0), ~\tau_j >t_0.
\end{equation}

We assume that for the initial function $\varphi$
the following hypothesis holds

(a4) $\varphi : (- \infty, t_0) \rightarrow {\bf R}$
is a Borel measurable bounded function.
\vspace{2 mm}

{\bf Definition.}
An absolutely continuous on each interval $[\tau_ j, \tau_{ j+1})$
function $x:[t_0, \infty)
\rightarrow {\bf R}$ is {\em a solution of the impulsive
problem (3),(4)}
if (3) is satisfied for almost all $t\in [0, \infty)$
and the equalities (4) hold.
\vspace{2 mm}

{\bf Definition.}
For each $s \geq 0$ the solution $X(t,s)$ of the problem
$$\dot{x}(t) + \sum_{k=1}^m A_k(t) x[h_k(t)]= 0,
\mbox{\rm~where~} t \geq s; ~x(\xi) = 0, ~\xi < s;$$
\begin{equation}
x(\tau_j)=B_j x(\tau_j -0), ~\tau_j >s, ~x(s)=1,
\end{equation}
{\em is a fundamental function} of the equation
(1),(2).

We assume $X(t,s) =0, ~0 \leq t < s$.

\newtheorem{uess}{Lemma}
\begin{uess}
$[12]$
Let (a1)-(a4) hold.
Then there exist one and only one solution
of the problem (3) with the initial condition $x(t_0)=\alpha_0$
and impulsive conditions
$$x(\tau_j) = B_j x(\tau_j)+\alpha_j$$
that can be presented in the form
$$ x(t) = X(t,t_0)x(t_0) +
\int_{t_0}^t X(t,s) f(s) ds -$$
\begin{equation}
- \sum_{k=1}^m \int_{t_0}^t X(t,s) A_k (s) \varphi [h_k(s)] ds
 +\sum_{\tau_j>t_0} X(t, \tau_j) \alpha_j,
\end{equation}
where $\varphi[h_k(s)]=0,$ if $h_k(s)>t_0$.
\end{uess}

\section{Non-oscillation Criteria for Impulsive Equations}

{\bf Definition.}
The equation (1),(2) {\em has a non-oscillating solution}
if there exist $t_0 >0,~ \varphi (t) $
satisfying (a4) such that for
$f \equiv 0$ the solution of (3),(4) is positive for
$t \geq t_0$.
Otherwise, all solutions of (1),(2) are said
to be {\em oscillating}.

In sequel we accept that the following hypothesis holds

(a5) delays are bounded: for every $s>0$
$$\mu_s = \min_k vrai\inf_{t>s} h_k(t) >-\infty$$
and there exists $s^{\prime}\geq s$ such that $ h_k(t)\geq s $
if $t\geq s^{\prime}.$

Denote for any $s$
\begin{equation}
\begin{array}{l}
A_k^{s}(t) =
\left\{
\begin{array}{ll}
A_k(t), & \mbox{~if~~~} t \geq s, \\
0,      & \mbox{~if~~~} t < s,
\end{array}
\right.    \\
h_k^{s}(t) =
\left\{
\begin{array}{ll}
h_k(t), & \mbox{~if~~~} t \geq s, \\
s,      & \mbox{~if~~~} t < s.
\end{array}
\right.  \end{array}
\end{equation}

The following theorem establishes non-oscillation criteria.

\newtheorem{guess}{Theorem}
\begin{guess}
Suppose (a1)-(a5) hold,
$A_k(t) \geq 0, ~k=1,\dots ,m,~ and
{}~B_j >0, ~j=1,2, \dots $.
Then the following hypotheses are equivalent

1) The equation (1),(2) has a non-oscillating
solution.

2) There exists $t_0 \geq 0$ such that
$X(t,s) > 0$, $t_0 \leq s < t < \infty$.

3) For a certain $t_1 \geq 0$ there exists
a non-negative integrable in each interval $[t_1,b]$
solution $u$ of an inequality
\begin{equation}
 u(t) \geq \sum_{k=1}^m A_k^{t_1}(t)
 \exp \left\{ \int_{h_k^{t_1} (t)}^t u(s) ds \right\}
\prod_{h_k^{t_1} (t) < \tau_j \leq t} B_j^{-1},
{}~t \geq t_1.
\end{equation}
Here and in sequel we assume that a product equals to unit
if number of factors is equal to zero.
\end{guess}

{\bf Proof.}
The scheme of the proof is
$1) \Longrightarrow 3) \Longrightarrow
2) \Longrightarrow 1). $
\vspace{2mm}

{\bf 1) $\Longrightarrow$ 3).}
Let  $x(t)$ be a positive solution
of (3),(4) ($f\equiv 0$).
By (a5) for a certain $t_1 \geq t_0 ~~h_k(t)>t_0,~~ t\geq t_1,
{}~k=1,\dots, m.$

Let us demonstrate that
$$u(t) = - \frac{d}{dt} \ln \left\{
\frac {x(t)}{x(t_1)} \prod_{t_1<\tau_j \leq t} B_j^{-1} \right\},
{}~~t \geq t_1.
$$
is a solution of (8).
To this end we integrate the latter equality
\begin{equation}
x(t) = x(t_1)\exp \left\{ - \int_{t_1}^t u(s) ds \right\}
\prod_{t_1 < \tau_j \leq t} B_j,~~ t \geq t_1.
\end{equation}

By setting $\varphi(t) = x(t)$ for $t<t_1$
one obtains that $x(t),~t \geq t_1$,
is a solution of (3),(4),
with the initial point $t=t_1$
and the initial function $\varphi (t)>0$.
We substitute (9) in (3) ($f\equiv 0 $):
$$ -u(t) \exp \left\{ - \int_{t_1}^t u(s) ds \right\}
\prod_{t_1<\tau_j \leq t} B_j + $$
$$\sum_{k \in N_1} A_k (t)
\exp \left\{ - \int_{t_1}^{h_k(t)} u(s) ds \right\}
\prod_{t_1<\tau_j \leq h_k(t)} B_j +
$$
\begin{equation}
\sum_{k \in N_2} A_k (t) \varphi[h_k (t)] = 0,
{}~t \geq t_1.
\end{equation}
Here $N_1 = \{k:~ h_k(t) \geq t_1 \},~~
N_2 = \{k:~ h_k(t) < t_1 \} $.

Using notations (7) the equality (10) can be rewritten in the form
$$ -u(t) \exp \left\{ - \int_{t_1}^t u(s) ds \right\}
\prod_{t_1< \tau_j \leq t} B_j +$$
$$
\sum_{k=1}^m A_k^{t_1} (t) \exp \left\{
- \int_{t_1}^{h_k^{t_1} (t)} u(s) ds \right\}
\prod_{t_1<\tau_j \leq h_k^{t_1} (t)} B_j +
$$
$$
\sum_{k \in N_2} A_k (t) \varphi[h_k (t)] =0,~
t \geq t_1. $$

Consequently,
$$
\left( u(t) - \sum_{k=1}^m A_k^{t_1} (t) \exp \left\{
 \int_{h_k^{t_1} (t)}^t u(s) ds \right\}
\prod_{h_k^{t_1} (t) < \tau_j \leq t} B_j^{-1} \right) \times
$$
$$
\times \exp \left\{
- \int_{t_1}^t u(s) ds \right\}
\prod_{t_1<\tau_j \leq t}  B_j =
\sum_{k \in N_2} A_k (t) \varphi[h_k (t)]  \geq  0,
$$
since $\varphi (t)$ is positive according
to our choice of the point $t_1$,
which implies 3).
\vspace{2mm}

{\bf 3) $\Longrightarrow $ 2).}
Consider (3),(4) with the initial function $\varphi
\equiv 0$ and initial value $x(t_1)=0$ in a segment
$[t_1,b]$:
\begin{equation}
\begin{array}{c}
\dot{x} (t) + \sum_{k=1}^m A_k (t) x[h_k(t)] =f(t),
{}~t \in [t_1,b]:~x(\xi) = 0, \xi < t_1,
\\
x(t_1)=0, ~x(\tau_j)=B_j x(\tau_j -0), ~\tau_j > t_1.
\end{array}
\end{equation}

Besides, we consider an ordinary impulsive differential equation
including the solution $u(t) \geq 0$ of (8):
$$ \dot{x} (t) +u(t) x(t) =z(t), ~t \in [t_1,b], $$
\begin{equation}
x(\tau_j)=B_j x(\tau_j-0), ~x(t_1)=0.
\end{equation}

The solution of (12) can be rewritten in the form [15]
\begin{equation}
x(t) = \int_{t_1}^t \exp \left\{ - \int_s^t
u(\xi) d\xi \right\}
\prod_{s<\tau_j\leq t} B_j z(s) ds.
\end{equation}

We seek for the solution of (11) of the form (13).
By substituting $x$ and $\dot{x}$ from (13) and (12)
into (11), we obtain
$$ z(t) - u(t) \int_{t_1}^t \exp \left\{
- \int_s^t u(\xi) d\xi \right\}
\prod_{s<\tau_j \leq t} B_j~ z(s) ds~ + $$
\begin{equation}
 \sum_{k=1}^m A_k^{t_1} (t) \int_{t_1}^{h_k^{t_1} (t)}
\exp \left\{
-\int_s^{h_k^{t_1}(t)}  u(\xi) d\xi \right\}
\prod_{s < \tau_j \leq h_k^{t_1}(t) } B_j~ z(s) ds = f(t).
\end{equation}

The equation (14) is of the type
\begin{equation}
z-Hz =f,
\end{equation}
where
$$(Hz)(t) = u(t) \int_{t_1}^t \exp \left\{
- \int_s^t u(\xi) d\xi \right\}
\prod_{s<\tau_j \leq t} B_j~ z(s) ds - $$
\begin{equation}
\sum_{k=1}^m A_k^{t_1} (t) \int_{t_1}^{h_k^{t_1} (t)}
\exp \left\{
-\int_s^{h_k^{t_1}(t)}  u(\xi) d\xi \right\}
\prod_{s < \tau_j \leq h_k^{t_1}(t) } B_j~ z(s) ds.
\end{equation}
It is well known [14] that the integral operator
$$ (Hz)(t) =
\int_{t_1}^b K(t,s)z(s) ds $$
acting in the space ${\bf L}_ {[t_1,b]}$
of functions integrable on $[t_1,b]$
is compact if
\begin{equation}
 | K(t,s) | \leq k(t), ~k \in {\bf L}_{ [t_1,b]}.
\end{equation}

For the operator $H$ defined by (16)
$$ |K(t,s)| \leq \sup_{s,t \in [t_1,b]}
\prod_{s < \tau_j \leq t} B_j
\left( u(t)
+ \sum_{k=1}^m |A_k (t) | \right) . $$

Thus the inequality (17) holds and
the operator $H: {\bf L}_{[t_1,b]}
\rightarrow {\bf L}_{[t_1,b]}$
is a compact Volterra integral operator.
Therefore [14] its spectral
radius is equal to zero.
Consequently the equation (15) for any $f \in {\bf L}_{[t_1,b]}$
has a single solution
\begin{equation}
z= (I-H)^{-1} f,
\end{equation}
where $I$ is the identity operator.

Let us show that $H$ is a positive operator.
The operator $H$ can be easily rewritten as a sum
$H= H_1 +H_2, $ where
$$ (H_1z)(t)=
\left( u(t) - \sum_{k=1}^m A_k^{t_1} (t)
\exp \left\{
\int_{h_k^{t_1}(t)}^t  u(s) ds \right\}
\prod_{h_k^{t_1} (t) < \tau_j \leq t}  B_j^{-1} \right)
\times $$ $$ \times \int_{t_1}^t \exp \left\{
- \int_s^t u(\xi) d\xi \right\}
\prod_{s<\tau_j \leq t} B_j ~z(s) ds , $$
$$ (H_2 z)(t)=
 \sum_{k=1}^m A_k^{t_1} (t) \int_{h_k^{t_1} (t)}^t
\exp \left\{
-\int_s^{h_k^{t_1}(t)}  u(\xi) d\xi \right\}
\prod_{s < \tau_j \leq h_k^{t_1} (t) }
 B_j~ z(s) ds . $$

The inequality (8) implies $H_1 \geq 0$.
So $H=H_1+H_2  \geq 0$.
Since the spectral radius of $H$ is equal to zero, then
$$ (I-H)^{-1} = I+H+H^2 + \dots \geq 0. $$
Thus if $f \geq 0,$ then
the solution $z$ of (15) is non-negative: $z \geq 0$.

The solution of (11) has the form (13),
where $z$ is the solution of (15).
Consequently, if in (11) $f \geq 0$,
then for the solution of (11) $x \geq 0$.

On the other hand, the solution of (11)
can be presented in the form (6)
\begin{equation}
x(t) = \int_{t_1}^t X(t,s)f(s)ds.
\end{equation}

As shown above, $f \geq 0$ implies $x \geq 0$.
Therefore the kernel of the integral operator
is non-negative, i.e. $X(t,s) \geq 0$ for
$t_1 \leq s \leq t <b$.
Since $b>t_1$ is chosen arbitrarily,
then $X(t,s) \geq 0$ for $t_1 \leq s <t < \infty$.

Let us prove that in fact the strict inequality
$X(t,s) >0$ holds.

Denote
$$x(t) = X(t,t_1) -\exp\left\{ - \int_{t_1}^t u(s) ds \right\}
\prod_{t_1<\tau_j \leq t} B_j . $$
Our purpose is to demonstrate $x(t)$ is non-negative.
The function $x(t)$ is a solution of (3),(4),
with $x(t_1)=0, \varphi \equiv 0$ and
$$ f(t) = u(t) \exp \left\{
- \int_{t_1}^t u(s) ds \right\}
\prod_{t_1<\tau_j \leq t} B_j - $$
$$
- \sum_{k=1}^m A_k^{t_1} (t)
\exp \left\{
-\int_{t_1}^{h_k^{t_1}(t)}  u(s) ds \right\}
\prod_{t_1 < \tau_j \leq h_k^{t_1}(t) } B_j = $$
$$ = \exp\left\{ - \int_{t_1}^t u(s) ds \right\}
\prod_{t_1<\tau_j \leq t} B_j \times $$ $$
\times \left( u(t) - \sum_{k=1}^m A_k^{t_1} (t)
\exp \left\{
\int_{h_k^{t_1}(t)}^t  u(s) ds \right\}
\prod_{h_k^{t_1}(t) < \tau_j \leq t } B_j^{-1} \right). $$

Thus (8) implies $f(t) \geq 0$.
Therefore in view of (6)
$$x(t) = \int_{t_1}^t X(t,s)f(s) ds \geq 0. $$
Consequently,
$$ X(t,t_1) \geq \exp\left\{ - \int_{t_1}^t u(s) ds \right\}
\prod_{t_1<\tau_j \leq t} B_j > 0. $$
For $s>t_1$ the inequality $X(t,s)>0$
can be proven similarly.
\vspace{2mm}

{\bf 2) $\Longrightarrow$ 1).}
Denote
$x(t)=X(t,t_0)$.
Then $x(t)$ is a positive solution
of (3),(4) ($f \equiv 0$) with
the initial function
$\varphi \equiv 0$.
The proof is complete.
\vspace{4mm}

Let us consider (1),(2) with coefficients
of an arbitrary sign.

Denote $a^+=\max\{a,0\},
a^- = \max \{-a,0 \}$.

\begin{guess}
Suppose (a1)-(a5) hold and $B_j>0 $.

Consider
three hypotheses:

1) The initial value problem (3),(4) with an initial point $t_0>0$
($f \equiv 0$) has a positive solution
that continuously extend the continuous initial function
$\varphi$.

2) $X(t,s) > 0, ~t_0 \leq s< t < \infty$.

3) There exists a non-negative integrable on each interval $[t_0,b]$
solution of an inequality
\begin{equation}
u(t) \geq \sum_{k=1}^m \left( A_k^{t_0} (t)\right)^+
\exp \left\{
\int_{h_k^{t_0}(t)}^t  u(s) ds \right\}
\prod_{h_k^{t_0}(t) < \tau_j \leq t } B_j^{-1} ,
{}~t \geq t_0.
\end{equation}

Then implications
$3) \Longrightarrow 2), 3) \Longrightarrow 1)$
are valid.
\end{guess}

{\bf Proof.}
The proof of ${\bf 3) \Longrightarrow 2)}$
coincides with the proof of $3) \Longrightarrow 2)$
in Theorem 1 up to the place where the
operator $H$ is presented as a sum of two terms.
Here
$$H=H_1+H_2+H_3,$$
where
$$ (H_1 z)(t)=
\left( u(t) - \sum_{k=1}^m \left( A_k^{t_0} (t)\right)^+
\exp \left\{
\int_{h_k^{t_0}(t)}^t  u(s) ds \right\}
\prod_{h_k^{t_0} (t) < \tau_j \leq t}  B_j^{-1} \right)
\times $$ $$ \times \int_{t_0}^t \exp \left\{
- \int_s^t u(\xi) d\xi \right\}
\prod_{s<\tau_j \leq t} B_j ~z(s) ds , $$
$$ (H_2 z)(t)=
 \sum_{k=1}^m \left( A_k^{t_0} (t)\right)^+
\int_{h_k^{t_0} (t)}^t
\exp \left\{
-\int_s^{h_k^{t_0}(t)}  u(\xi) d\xi \right\}
\prod_{s < \tau_j \leq h_k^{t_0} (t) }
 B_j~ z(s) ds , $$
$$ (H_3 z)(t)=
 \sum_{k=1}^m \left( A_k^{t_0} (t)\right)^-
\int_{t_0}^{h_k^{t_0} (t)}
\exp \left\{
-\int_s^{h_k^{t_0}(t)}  u(\xi) d\xi \right\}
\prod_{s < \tau_j \leq h_k^{t_0} (t) }
 B_j~ z(s) ds . $$
Again, like in Theorem 1, $H_1 \geq 0,~H_2 \geq 0,
{}~H_3 \geq 0$, which implies
$H=H_1+H_2+H_3 \geq 0$.
The end of the proof completely repeats
the corresponding one of Theorem 1.

${\bf 3) \Longrightarrow 1)}$.
Let us consider the problem (3),(4) .
Let $ \mu _{t_0} $ be chosen as in the hypothesis (a5).
We extend to the interval
$[\mu_{t_0}, t_0 )$ the coefficients $A_k(t)$ by zero
and the delays $h_k(t)$ such that $h_k(t) \leq t$.
Let $u(t)$ be a non-negative function
satisfying (20).
We extend it by zero to $[\mu_{t_0}, t_0 )$.
Then $u(t)$ is a solution of (20), where $t_0$
is changed by $\mu_{t_0}$.

Consider a corresponding extension of (3),(4)
to the interval $[\mu_{t_0}, \infty )$.
As proven above, 3) $\Longrightarrow$ 2),
therefore $X(t,s) >0$
for $\mu_{t_0} \leq s < t < \infty$.
Assuming
$$\varphi(t) = X(t, \mu_{t_0} ) \mbox{~for~~~}
\mu_{t_0} \leq t < t_0 \mbox{~~~~and} $$
$$ x(t) = X(t, \mu_{t_0} )  \mbox{~for~~~} t \geq t_0 ,$$
we obtain that $x(t)$ is a positive solution of (3),(4)
($f \equiv 0$), with an initial point $t_0$,
that continuously extends
the continuous initial function $\varphi$.
This completes the proof of the theorem.
\vspace{3mm}

Now we proceed to explicit non-oscillation results.

Denote
$$\underline{h}^{t_0}(t)=\min _k h_k^{t_0}(t),$$
where $h_k^s (t)$ is defined by (7).

\begin{guess}
Suppose (a1)-(a5) hold,
$B_j>0$ and at least one of the following
three hypotheses hold:

{}~~1) $A_k(t) \leq 0, ~t \geq t_0.$
\begin{equation}
2)~~~~~~vrai\sup_{t>t_0}
\sum_{k=1}^m
\int_{\underline{h}^{t_0} (t)}^t
\left( A_k^{t_0} (s)\right)^+
\prod_{h_k^{t_0} (s) < \tau_j \leq s }
 B_j^{-1}~ z(s) ds \leq 1/e.
\end{equation}
\begin{equation}
3)~
\sum_{k=1}^m
\int_{\underline{h}^{t_0} (t)}^t
\left( A_k^{t_0} (s)\right)^+  ds
\leq 1/e \left( 1+
{\sum_{h_k^{t_0} (t) < \tau_j \leq t ,~ B_j < 1 }}
\!\!\!\! \ln B_j \right),
{}~~~t \geq t_0.
\end{equation}

Then the fundamental matrix $X(t,s)$
is positive for $t_0 \leq s < t < \infty$
and there exists a positive solution
of (3),(4) ($f \equiv 0$)
continuously extending a continuous initial function
$\varphi$.
\end{guess}

{\bf Proof.}
Obviously 1) is a special case of 2).
Let us prove the theorem assuming (21) holds.
To this end we will demonstrate that
a function
$$
u(t) = e \sum_{k=1}^m \left( A_k^{t_0} (t)\right)^+
\prod_{h_k^{t_0} (t) < \tau_j \leq t }
 B_j^{-1}
 $$
is a non-negative solution of the inequality (20).
By substituting $u$ in (20) one obtains
$$
 e \sum_{k=1}^m \left( A_k^{t_0} (t)\right)^+
\prod_{h_k^{t_0} (t) < \tau_j \leq t }
 B_j ^{-1} \geq
\sum_{k=1}^m \left( A_k^{t_0} (t)\right)^+
\times $$
$$
\times \exp \left\{ e
\int_{h_k^{t_0} (t)}^t
\sum_{i=1}^m
\left( A_i^{t_0} (s)\right)^+
\prod_{h_i^{t_0} (s) < \tau_j \leq s }
 B_j^{-1}~  ds
 \right\}
\prod_{h_k^{t_0} (t) < \tau_j \leq t }
 B_j ^{-1} . $$

This unequality can be deduced from the following one
$$
 e \sum_{k=1}^m \left( A_k^{t_0} (t)\right)^+
\prod_{h_k^{t_0} (t) < \tau_j \leq t }
 B_j ^{-1} \geq  $$ $$
\exp \left\{ e
\int_{\underline{h}^{t_0} (t)}^t
\sum_{k=1}^m
\left( A_k^{t_0} (s)\right)^+
\prod_{h_k^{t_0} (s) < \tau_j \leq s }
 B_j^{-1}~  ds
 \right\} \times $$ $$\times
 \sum_{k=1}^m \left( A_k^{t_0} (t)\right)^+
\prod_{h_k^{t_0} (t) < \tau_j \leq t }
 B_j ^{-1} . $$
After dividing this inequality by its left-hand side
and logarithmizing it we obtain
$$\sum_{k=1}^m
\int_{\underline h^{t_0} (t)}^t
\left( A_k^{t_0} (s)\right)^+
\prod_{h_k^{t_0} (s) < \tau_j \leq s }
 B_j^{-1}~ z(s) ds \leq 1/e, $$
which obviously results from (21).
\vspace{2mm}

Let 3) hold.
We will prove that
$$
u(t) = e \sum_{k=1}^m \left( A_k^{t_0} (t)\right)^+
$$
is a solution of the inequality (20) which after
substituting takes form
$$
 e \sum_{k=1}^m \left( A_k^{t_0} (t)\right)^+
 \geq  $$ $$
\sum_{k=1}^m \left( A_k^{t_0} (t)\right)^+
\exp \left\{ e
\int_{h_k^{t_0} (t)}^t
\sum_{i=1}^m
\left( A_i^{t_0} (s) \right)^+ ds
 \right\}
\prod_{h_k^{t_0} (t) < \tau_j \leq t }
 B_j ^{-1} . $$
This inequality can be deduced from
$$
 e \sum_{k=1}^m \left( A_k^{t_0} (t)\right)^+
 \geq  $$ $$
\sum_{k=1}^m \left( A_k^{t_0} (t)\right)^+
\exp \left\{ e
\int_{\underline {h}^{t_0} (t)}^t
\sum_{k=1}^m
\left( A_k^{t_0} (s)\right)^+ ds
 \right\}
\prod_{\underline {h}^{t_0} (t) < \tau_j \leq t }
 B_j ^{-1} , $$
where the product contains only factors for which $B_j<1.$
The latter inequality after dividing
by the left-hand side and logarithmizing
coincides with (22).
This completes the proof of the theorem.
\vspace{2mm}

Let us compare oscillation properties of (1),(2) and
an impulsive equation
$$ \dot{x} (t) + \sum_{k=1}^m \tilde{A}_k (t)
x(\tilde{h}_k (t)) = f(t), ~t \in [0, \infty), $$
\begin{equation}
x(\tau_j)=\tilde{B}_j x(\tau_j-0).
\end{equation}

\begin{guess}
Let the hypotheses (a1)-(a5) hold for the equations (1),(2)
and (23), $\tilde {A}_k(t)\geq 0, B_j >0$.
Suppose that any (therefore, all) of the
hypotheses 1)-3) of Theorem 1 holds for (1),(2).

Then if $A_k (t) \geq \tilde{A}_k (t), ~ B_j \leq \tilde{B}_j$
and at least one of the hypotheses

1) $h_k (t) \leq \tilde{h}_k (t), ~\tilde{B}_j \leq 1,
{}~j=1,2, \dots ;$

2) $h_k (t) = \tilde{h}_k(t)$,

holds then for the equation (23) the assertions 1)-3)
of Theorem 1 are valid.

\end{guess}

{\bf Proof.}
By the hypothesis of the theorem
there exists a non-negative function $u(t)$
satisfying (7).
Besides, for any non-negative function $u$
under the hypotheses of the theorem
the inequality
$$
 \sum_{k=1}^m A_k (t)
\exp\left\{ \int_{h_k (t)}^t
u (s) ds \right\}
\prod_{h_k(t) < \tau_j \leq t} B_j^{-1} \geq $$ $$ \geq
\sum_{k=1}^m \tilde{A}_k (t)
\exp\left\{ \int_{\tilde{h}_k (t)}^t
u(s) ds \right\}
\prod_{\tilde{h}_k(t) < \tau_j \leq t} \tilde{B}_j^{-1}
$$ holds.
Consequently if $u$ is a solution of the inequality (7)
then $u$ is a solution of this inequality,
where $A_k, h_k, B_j$ are changed by $\tilde{A}_k,
\tilde{h}_k, \tilde{B}_j$.
Then by Theorem 1 the other
assertions of this theorem also hold .
\vspace{2 mm}

{\bf Corollary 1.}
Suppose the hypotheses (a1)-(a5) hold for (1),(2)
and $B_j >0$.
Besides, let $0 \leq A_k (t) \leq A_k, ~t-h_k(t) \leq h_k,
{}~ B_j \leq 1$.

If there exists a non-oscillating solution
of the equation with constant coefficients and delays
$$
\dot{x} (t) + \sum_{k=1}^m A_k x(t-h_k)=f(t),
{}~t \in [0, \infty), $$
$$
{}~x(\tau_j) = B_j x(\tau_j -0), $$
then there exists a non-oscillating solution
of the equation (1),(2).
\vspace{2 mm}

{\bf Corollary 2.}
Let (a1)-(a5) hold and $A_k (t) \geq 0.$
If there exists a non-oscillation solution
of the equation (1) without impulses
and $B_j \geq 1$, then there exists
a non-oscillating solution of
the impulsive equation (1),(2).
\vspace{2mm}

\section{Oscillation Properties of Impulsive and
Non-impulsive Equations}

Consider a non-impulsive differential equation
\begin{equation}
\dot{x} (t) + \sum_{k=1}^m a_k(t) x[h_k(t)] = f(t),
{}~ t \geq 0.
\end{equation}
Denote by $x(t,s)$ the fundamental function of the equation (24).
After substituting $B_j \equiv 1$
Theorems 1 and 2 immediately yield the following results.

\begin{guess}
Suppose (a2)-(a5) hold for (24) and $a_k(t) \geq 0, k=1,2, \dots$.
Then the following hypotheses are equivalent:

1) The equation (24) has a non-oscillating
solution $(f\equiv 0)$.

2) There exists $t_0 \geq 0$ such that
$x(t,s) >0$ for $t_0 \leq s < t < \infty$.

3) For a certain $t_1 \geq 0$ there exists
a non-negative integrable on each interval $[t_1,b]$
solution $u$ of the inequality
\begin{equation}
u(t) \geq \sum_{k=1}^m a_k^{t_1}(t)
\exp \left\{ \int_{h_k^{t_1}(t)}^t u(\xi) d\xi \right\},
{}~t \geq t_1.
\end{equation}
\end{guess}

\begin{guess}
Suppose (a2)-(a5) hold for (24).
Consider three hypotheses:

1) The initial value problem for (24)
($f \equiv 0$) with an initial point $t_0 \geq 0$
has a positive solution that is a continuous expansion
of a continuous initial function $\varphi$;

2) $x(t,s)>0, ~t_0 \leq s < t < \infty$;

3) There exists a non-negative integrable on each interval
$[t_0,b]$ solution $u$ of the inequality
$$u(t) \geq \sum_{k=1}^m \left(
a_k^{t_0} (t) \right)^+ \exp \left\{
\int_{h_k^{t_0} (t)}^{t} u(s) ds \right\},
{}~t \geq t_0 $$.

Then implications  $ 3) \Rightarrow 2) ,3\Rightarrow 1)$
are valid.
\end{guess}

Corollary of Theorem 3 for the equation (24)
coincides with the known non-oscillation result
for equations without impulses [1,2,4].

In this paper we present a fundamental result
that enables to reduce the oscillation problem for (1),(2)
to the oscillation problem for an equation without impulses.
To this end consider an auxiliary equation
\begin{equation}
\dot{x} (t) + \sum_{k=1}^m A_k(t)
\prod_{h_k^0 (t) < \tau_j \leq t} B_j^{-1}
x[h_k(t)]= 0, ~t \in [0, \infty),
\end{equation}
$$ \mbox{where~~~~~~}
h_k^0 (t) =
\left\{
\begin{array}{ll}
h_k(t),   &  \mbox{~if~~} t \geq 0, \\
0,        &  \mbox{~if~~} t <0.
\end{array}
\right.
$$

Denote by $Y(t,s)$ a fundamental function
of the equation (26).

\begin{guess}
Suppose (a1)-(a5) hold, $A_k \geq 0, ~B_j >0$.

Then

1) There exists $t_0>0$, such that $X(t,s) > 0, ~t_0 \leq s < t < \infty$
$\Longleftrightarrow$ there exists $t_1>0$, such that $Y(t,s)> 0,
{}~t_1 \leq s<t< \infty$.

2) All solutions of (1),(2) ($f \equiv 0$)
are oscillating $\Longleftrightarrow $
all solutions of (26)  are oscillating.

3) There exists a non-oscillating solution
of (1),(2) ($f \equiv 0$)
$\Longleftrightarrow$ there exists
a non-oscillating solution
of (26).
\end{guess}

{\bf Proof.}
1). Let $X(t,s) > 0, ~t_0 \leq s < t < \infty$.
Then by Theorem 1 there exists a solution
of the inequality (7)
for $t \geq t_1$.
This inequality coincides with (25) under
$$a_k (t) = A_k (t)
\prod_{h_k^0 (t) < \tau_j \leq t}  B_j^{-1} . $$
Therefore by Theorem 5 $Y(t,s) >0,~ t_1\leq s <t <\infty$.
The converse can be proven similarly.

2). Suppose all solutions of (1),(2)
($f \equiv 0$) are oscillating and (26)
has a positive solution, beginning
with a certain $t_0$.
Then by Theorem 5 $Y(t,s) >0$
for $t_1 \leq s <t < \infty$.
Then, as proven in 1), $X(t,s) >0$
for $t_2 \leq s < t < \infty$.
Consequently, by Theorem 1  the equation (1),(2)
has a non-oscillating solution,
which contradicts to the hypothesis.
The converse is proven similarly.

Besides, 2) implies 3), which completes the proof.

By applying Theorem 7 and known
oscillation (non-oscillation) results
on equations without impulses,
one obtains oscillation results for impulsive equations.
As an example we present the following statement.

$$ \mbox{Denote~~~} \underline{h} (t) = \min_k h_k(t), ~
\bar{h} (t) = \max_k h_k(t). $$

\begin{guess}
Let (a1)-(a5) hold for (1),(2), $A_k(t)\geq 0$ and $B_j >0$.
Then if at least one of the following
inequalities holds
$$ 1) ~ \lim_{t \rightarrow \infty}\inf
\int_{\underline{h} (t)}^t \sum_{k=1}^m A_k (s)
\prod_{h_k(s) < \tau_j \leq s} B_j^{-1} ds > 1/e,
$$
$$
2) ~ \lim_{t \rightarrow \infty}\sup
\int_{\bar {h} (t)}^t
\sum_{k=1}^m A_k (s)
\prod_{h_k(s) < \tau_j \leq s} B_j^{-1} ds > 1,
$$
then all the solutions of (1),(2) are oscillating.
\end{guess}

This statement is obtained by applying Theorem 7
and oscillation results for equations without
impulses from the monographs [1,2,4].

\end{document}